\documentclass[reqno,12pt]{amsart}

\usepackage[ps,dvips,all]{xy}

\usepackage{latexsym,amssymb,amsmath,amsfonts,epsfig,textcomp}
\usepackage{pstricks,pst-coil}

\def\sqr#1#2{{\vcenter{\hrule height.#2pt
        \hbox{\vrule width.#2pt height#1pt \kern#1pt
                \vrule width.#2pt}
        \hrule height.#2pt}}}

\newtheorem{theorem}{Theorem}[section]
\newtheorem{lemma}[theorem]{Lemma}
\newtheorem{proposition}[theorem]{Proposition}
\newtheorem{corollary}[theorem]{Corollary}

\theoremstyle{definition}
\newtheorem{definition}[theorem]{Definition} % \theoremstyle{remark}
\newtheorem{remark}[theorem]{Remark}
\newtheorem{example}[theorem]{Example}
\newtheorem{assumption}[theorem]{Assumption}

\newcommand {\ZZ}{\mathbb{Z}}
\newcommand {\NN}{\mathbb{N}}
\newcommand {\PP}{\mathbb{P}}
\newcommand {\RR}{\mathbb{R}}

\newcommand{\mif}{\mbox{if} ~}
\newcommand{\s}{\; | \;}

\newcommand{\ov}[1]{\overline{#1}}

\newcommand{\lm}{_{\lambda -\mu}}

\DeclareMathOperator{\reg}{reg}
\DeclareMathOperator{\htt}{ht}

\DeclareMathOperator{\depth}{depth}

\begin{document}

\title[ Specializations of Ferrers ideals]
{ Specializations of Ferrers ideals}

\author[A.\ Corso and U.\ Nagel]
{Alberto Corso \and Uwe Nagel$^{*\, \dagger}$}

%\thanks{AMS 2000 {\em Mathematics Subject Classification}.
%Primary:  05A15, 13D02, 13D40, 14M25; Secondary: 05C75, 13C40, 13H10,  14M12,
% 52B05. ?????}

\thanks{$^*$ This author gratefully acknowledges partial support
from the NSA under grant H98230-07-1-0065}
\thanks{$^\dagger$ Current address: Institute for Mathematics \& its Applications,
University of Minnesota, Minneapolis, MN 554555, USA}

\address{Department of Mathematics, University of Kentucky, Lexington, KY 40506, USA}
\email{corso@ms.uky.edu}
\email{uwenagel@ms.uky.edu}

\begin{abstract}
We introduce a specialization technique in order to study monomial
ideals that are generated in degree two by using our earlier results
about Ferrers ideals. It allows us to describe explicitly a cellular
minimal free resolution of various ideals including any strongly
stable and any squarefree strongly stable ideal whose minimal
generators have degree two. In particular, this shows that threshold
graphs can be obtained as specializations of Ferrers graphs, which
explains their similar properties.
\end{abstract}

%\pagestyle{empty}
% \subjclass{Primary ???; Secondary ???}

\maketitle

\vspace{0.1in}

\section{Introduction}
One of the starting points of this note has been the observation
that two very common classes of graphs, namely Ferrers graphs and
threshold graphs, have similar properties (see, e.g., \cite{MP}). This is remarkable as
Ferrers graphs are particular bipartite graphs on vertex sets
$\{x_1,\ldots,x_n\}$ and $\{y_1,\ldots,y_m\}$, whereas threshold
graphs are typically not bipartite. One of the goals of this note is
to show that the similarity between these graphs extends to
algebraic properties of their edge ideals and that it has a natural
interpretation. In fact, in \cite{CN1} we have described a cellular
minimal free resolution of Ferrers ideals, the edge ideals of
Ferrers graphs. The polyhedral cell complex that governs this
cellular resolution has a very nice geometric description as a
certain subcomplex of the face complex of the product of two
simplices. This allows us to compute various invariants of Ferrers
ideals as, for example, their $\ZZ$-graded Betti numbers and their
height. The main idea of this note is that this information can be
used to obtain insight about graphs that are often not bipartite by
using a `specialization' process (see Section~\ref{sec-spec}).
Roughly speaking, specializing simply means to identify each
$y$-vertex with an $x$-vertex. Extending this specialization to the polyhedral
cell complex that resolves the Ferrers ideal provides, under
suitable hypotheses, a cellular minimal free resolution of the
specialized Ferrers ideal that is not necessarily a squarefree
monomial ideal. After some preliminaries this program is carried out
in Section~\ref{sec-spec}. In Section~\ref{sec-thres} we discuss the
class of ideals and graphs that can be described by using
specializations of Ferrers ideals. In particular, we show that all
threshold graphs can be obtained as specializations of Ferrers
graphs. Furthermore, every strongly stable ideal that is generated
in degree two can be obtained as such a specialization.

\smallskip

Horwitz shows in \cite{H} that each squarefree monomial ideal $I$
that has a $2$-linear free resolution admits a cellular minimal free
resolution that is given by a regular cell complex,
provided the graph that $I$ corresponds to does not
contain a certain subgraph $G'$. Though we consider only a subset of
the monomial ideals with regularity two, our results for these are
more explicit. In particular, we give a geometric description of the
underlying polyhedral cell complexes. We also show that our results apply to
the exceptional graph $G'$ (see Example~\ref{ex-exception}).

%%%%%%%%%%%%%%%%%%%%%%%%%%%%%%%%%%%%%%%%%%%%%%%%%%%%%%%%%%%%%%%%%%%%%%%%%%%%%%%%%%%%%%%%%

\section{Preliminaries} \label{sec-prel}

A {\it Ferrers graph} is a bipartite graph on two distinct vertex
sets ${\bf X}=\{ x_1, \ldots, x_n\}$ and ${\bf Y} =\{ y_1, \ldots,
y_m \}$ such that if $(x_i, y_j)$ is an edge of $G$, then so is
$(x_p, y_q)$ for $1 \leq p \leq i$ and $1 \leq q \leq j$. In
addition, $(x_1,y_m)$ and $(x_n,y_1)$ are required to be edges of
$G$. For any Ferrers graph $G$ there is an associated sequence of
non-negative integers $\lambda = (\lambda_1, \lambda_2, \ldots,
\lambda_n)$, where $\lambda_i$ is the degree of the vertex $x_i$.
Notice that the defining properties of a Ferrers graph imply that
$\lambda_1=m \geq \lambda_2 \geq \cdots \geq \lambda_n \geq 1$; thus
$\lambda$ is a {\it partition}. Alternatively, we can associate to a
Ferrers graph a diagram ${\mathbf T}_{\lambda}$, dubbed {\it Ferrers
tableau}, consisting of an array of $n$ rows of cells with
$\lambda_i$ adjacent cells, left justified, in the $i$-th row. A
{\it Ferrers ideals} is the edge ideal associated with a Ferrers
graph. See \cite{CN1} and \cite{vila}  for additional details.

\medskip

In this note we study ideals that are closely related to Ferrers
ideals. In order to explicitly describe their minimal free
resolutions, we use the theory of cellular resolutions and
polyhedral cell complexes as developed in \cite{BPS} and \cite{BS}.
We briefly recall some basic notions.  However we refer to \cite{BS}
$($or \cite{MS}$)$ for a more detailed introduction.

A {\it polyhedral cell complex}
$X$ is a finite collection of convex polytopes $($in some ${\mathbb
R}^N)$ called faces $($or cells$)$ of $X$ such that:
\begin{enumerate}
\item
if $P \in X$ and $F$ is a face of $P$, then $F \in X$;

\item
if $P, Q \in X$ then $P \cap Q$ is a face of both $P$ and $Q$.
\end{enumerate}
Let $F_k(X)$ be the set of $k$-dimensional faces. Each cell complex
admits an {\it incidence function} $\varepsilon$ on $X$, where
$\varepsilon(Q,P) \in \{ 1, -1 \}$ if $Q$ is a facet of $P \in X$.
$X$ is called a {\it labeled} cell complex if each vertex $i$ has a
vector ${\mathbf a}_i \in {\mathbb N}^N$ $($or the monomial
${\mathbf z}^{{\mathbf a}_i}$, where ${\mathbf z}^{{\mathbf a}_i}$
denotes a monomial in the variables $z_1, \ldots, z_N)$ as label.
The label of an arbitrary face $Q$ of $X$ is the exponent ${\mathbf
a}_Q$, where ${\mathbf z}^{{\mathbf a}_Q} := {\rm lcm}\, ({\mathbf
z}^{{\mathbf a}_i} \,|\, i \in Q)$. Each labeled cell complex
determines a complex of free $R$-modules, where $R$ is the
polynomial ring $K[z_1,\ldots,z_N]$. The cellular complex ${\mathcal
F}_X$ supported on $X$ is the complex of free $\ZZ^N$-graded
$R$-modules
\[
{\mathcal F}_X\colon \quad 0 \rightarrow S^{F_d(X)}
\stackrel{\partial_d}{\longrightarrow} S^{F_{d-1}(X)}
\stackrel{\partial_{d-1}}{\longrightarrow} \cdots
\stackrel{\partial_2}{\longrightarrow} S^{F_1(X)}
\stackrel{\partial_1}{\longrightarrow} S^{F_0(X)}
\stackrel{\partial_0}{\longrightarrow} S \rightarrow 0,
\]
where $d= \dim X$ and $S^{F_k(X)}:= \displaystyle \bigoplus_{P \in
F_k(X)} R[-{\mathbf a}_P]$. The map $\partial_k$ is defined by
\[
\partial_k(e_P) := \sum_{{Q {\rm \ facet \ of \ } P}}
\varepsilon(P,Q) \cdot {\mathbf z}^{{\mathbf a}_P -{\mathbf a}_Q}
\cdot e_Q,
\]
where $\{ e_P \,|\, P \in F_k(X) \}$ is a basis of $S^{F_k(X)}$ and
$e_{\emptyset} :=1 $. If ${\mathcal F}_X$ is acyclic, then it
provides a free $\ZZ^N$-graded resolution of the image $I$ of
$\partial_0$, that is the ideal generated by the labels of the
vertices of $X$. In  this case, ${\mathcal F}_X$ is called a {\it
cellular resolution} of $I$.

\begin{example} \label{ex-compl-bip}
Consider the ideal $ I := (x_1, \ldots, x_n)(y_1, \ldots, y_m)
\subset R = K[x_0, \ldots, x_n, y_1,$ $\ldots, y_m]$. Let $X_{n, m}$
be the face complex of the polytope $\Delta_{n-1} \times
\Delta_{m-1}$ obtained by taking the cartesian product of the
$(n-1)$-simplex $\Delta_{n-1}$ and the $(m-1)$-simplex
$\Delta_{m-1}$. Labeling the vertices of $\Delta_{n-1}$ by
$x_1,\ldots,x_n$ and the ones of $\Delta_{m-1}$ by $y_1,\ldots,y_m$,
the vertices of the cell complex $X_{n,m}$ are naturally labeled by
the monomials $x_iy_j$ with $1 \leq i \leq n$ and $1 \leq j \leq m$.
This turns $X_{n, m}$ into a labeled polyhedral cell complex. The
picture below illustrates the case where $n=2$ and $m=3$. It is
shown in \cite{CN1} that the complex ${\mathcal F}_{X_{n,m}}$ is a
minimal free resolution of $I$.

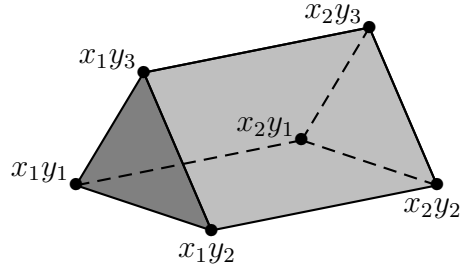
\begin{center}
\begin{pspicture}(5,0)(11,4)
\psset{xunit=.6cm, yunit=.6cm}

\pscustom[linestyle=none,fillstyle=solid,fillcolor=gray]{
\psline(9,2)(10.5,4.5) \psline[liftpen=0](10.5,4.5)(12,1) }

\pscustom[fillstyle=solid,fillcolor=lightgray]{
\psline(12,1)(10.5,4.5) \psline[liftpen=0](15.5,5.5)(17,2) }

\rput(9,2){$\bullet$} \rput(12,1){$\bullet$}
\rput(10.5,4.5){$\bullet$} \rput(14,3){$\bullet$}
\rput(17,2){$\bullet$} \rput(15.5,5.5){$\bullet$}

\rput(8.2,2.3){$x_1y_1$} \rput(11.9,0.5){$x_1y_2$}
\rput(9.7,4.8){$x_1y_3$} \rput(13.2,3.3){$x_2y_1$}
\rput(16.9,1.5){$x_2y_2$} \rput(14.7,5.8){$x_2y_3$}

\psline[linestyle=solid](9,2)(12,1)
\psline[linestyle=solid](9,2)(10.5,4.5)
\psline[linestyle=solid](12,1)(10.5,4.5)
\psline[linestyle=solid](17,2)(15.5,5.5)
\psline[linestyle=solid](12,1)(17,2)
\psline[linestyle=solid](10.5,4.5)(15.5,5.5)
\psline[linestyle=dashed](14,3)(17,2)
\psline[linestyle=dashed](14,3)(15.5,5.5)
\psline[linestyle=dashed](9,2)(14,3)
\end{pspicture}
\end{center}

\end{example}

The acyclicity of ${\mathcal F}_X$  is merely determined by the
geometry of the polyhedral cell complex $X$. Recall that $X$ is
called {\em acyclic} if it is either empty or has zero reduced
homology. Moreover, consider the partial order on $\NN_0^N$ defined
by ${\mathbf a} \preccurlyeq {\mathbf c}$ if ${\mathbf c} - {\mathbf
a} \in \NN_0^N$. For any ${\mathbf c} \in \ZZ^N$, we define the
subcomplex $X_{\preccurlyeq {\mathbf c}}$ of $X$ as the labeled
complex that consists of the faces of $X$ whose labeling monomials
$z^{\mathbf a}$ satisfy ${\mathbf a} \preccurlyeq {\mathbf c}$.

Bayer and Sturmfels \cite{BS}, Proposition 1.2, have established the
following criterion that we will use in the following section.

\begin{lemma} \label{lem-acyc-crit}
The complex ${\mathcal F}_X$ is a cellular resolution  if and only
if, for each ${\mathbf c} \in \NN_0^N$, the complex $X_{\preccurlyeq
{\mathbf c}}$ is acyclic over the field $K$.
\end{lemma}

%%%%%%%%%%%%%%%%%%%%%%%%%%%%%%%%%%%%%%%%%%%%%%%%%%%%%%%%

\section{Specializations} \label{sec-spec}

There are relatively few monomial ideals for which the minimal free
resolution is explicitly known. These include the edge ideals of
bipartite graphs that are 2-regular. Up to isomorphisms these are
exactly the Ferrers ideals whose minimal free resolutions have been
described in \cite{CN1}. Here we want to show that this information
can be used to obtain the minimal free resolution of other monomial
ideals by a process that we call specialization. This resolution
will be again cellular.

\begin{definition} \label{def-spec}
Let $I$ be a monomial ideal contained in $R = K[x_0, \ldots, x_n,
y_1, \ldots, y_m]$. Let $\sigma\colon \{y_1,\ldots,y_m\}
\longrightarrow \{x_1,\ldots,x_k\}$ be any map, where $k = \max \{m, n\}$
and  $x_{n+1}, \ldots, x_k$
are (possibly) additional variables. By abuse of notation
we use the same symbol to denote the substitution homomorphism
$\sigma\colon R \longrightarrow S$, where $S := K[x_1,\ldots,x_k]$,
given by $x_i \mapsto x_i$ and $y_i \mapsto \sigma (x_i)$. We call
$\sigma$ the {\it specialization map} and the monomial ideal
$\overline{I} := \sigma (I) \subset S$ the {\it specialization of
$I$}.
\end{definition}

In general, the ideals $I$ and $\overline{I}$ have quite different
properties.

\begin{example} \label{ex-bad-spec}
(i) Let $\lambda := (2, 2)$ and consider the specialization $\sigma$
defined by $y_i \mapsto x_i$. Then the Ferrers ideal $I_{\lambda}$
has 4 minimal generators while $\overline{I}_{\lambda}$ has only 3
minimal generators.

(ii) Consider the ideal $I = (x_1 y_1, x_1 y_3, x_2 y_1)$ and the
specialization $\sigma (y_i) = x_i$. Then  $I$ and its
specialization $\overline{I}$ have the same number of minimal
generators, but $I$ has height two whereas $\overline{I}$ has height
one. However, if we use the specialization defined by $y_i \mapsto
x_{4-i}$, then $I$ and the specialized ideal have the same
$\ZZ$-graded Betti numbers.
\end{example}

These examples illustrate that we need some assumptions and a
careful choice of the specialization in order to study the
specialized ideal by means of the original one. Throughout the
remainder of this note we make the following

\begin{assumption}
Let us assume that $m \geq n$ and that $\sigma\colon R \longrightarrow
S :=K[x_1, \ldots, x_m]$ is the specialization map
defined by
\[
\sigma (y_i) = x_i.
\]
\end{assumption}

In order to increase the range of graphs obtained as a
specialization of Ferrers graphs, we introduce some notation for
ideals that are isomorphic to Ferrers ideals.

\begin{definition} \label{def-genFer}
Let $\lambda = (\lambda_1,\ldots,\lambda_n)$ be a partition and let
$\mu = (\mu_1,\ldots,\mu_n) \in \ZZ^n$ be a vector such that
\[
0 \leq \mu_1 \leq \cdots \leq \mu_n < \lambda_n.
\]
Then we define the ideal
\[
I_{\lambda - \mu} := (x_i y_j \s 1 \leq i \leq n, \mu_i < j \leq
\lambda_i) \subset I_{\lambda}
\]
and call it a {\it generalized Ferrers ideal}.
\end{definition}

Note that the assumption $\mu_n < \lambda_n$ is essentially not a
restriction. It just ensures that the variable $x_n$ divides one of
the minimal generators of $I_{\lambda - \mu}$.

As in the case of Ferrers ideals, generalized Ferrers ideals
correspond to a shape ${\mathbf T}_{\lambda - \mu}$ that is obtained
from the Ferrers diagram ${\mathbf T}_{\lambda}$ by removing the
first $\mu_i$ boxes in row $i$ beginning on the left-hand side. We
use the notation $\lambda - \mu$ in order to distinguish it from the
common notation for skew shapes. Two examples are illustrated below:

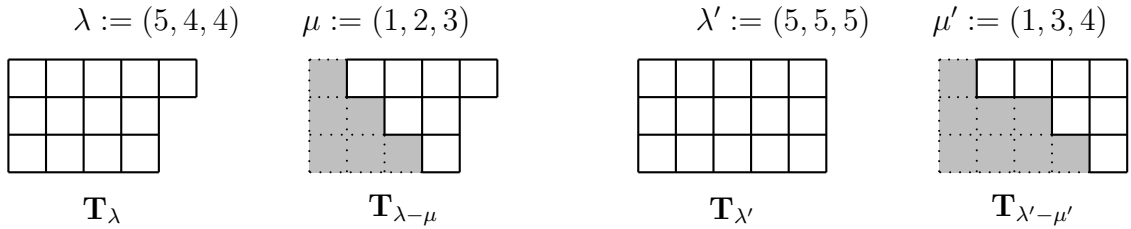
\begin{center}
\begin{pspicture}(0,0)(7,4)
\psset{xunit=.5cm, yunit=.5cm}

\rput(7,6){$\lambda:=(5,4,4) \quad \quad \mu:=(1,2,3)$}

\psline[linestyle=solid](0,5)(5,5)
\psline[linestyle=solid](0,4)(5,4)
\psline[linestyle=solid](0,3)(4,3)
\psline[linestyle=solid](0,2)(4,2)

\psline[linestyle=solid](0,2)(0,5)
\psline[linestyle=solid](1,2)(1,5)
\psline[linestyle=solid](2,2)(2,5)
\psline[linestyle=solid](3,2)(3,5)
\psline[linestyle=solid](4,2)(4,5)
\psline[linestyle=solid](5,4)(5,5)

\rput(2.5,1){${\bf T}_{\lambda}$}

\pscustom[linestyle=none,fillstyle=solid,fillcolor=lightgray]{
\psline(8,5)(9,5) \psline[liftpen=0](9,4)(8,4) }

\pscustom[linestyle=none,fillstyle=solid,fillcolor=lightgray]{
\psline(8,4)(10,4) \psline[liftpen=0](10,3)(8,3) }

\pscustom[linestyle=none,fillstyle=solid,fillcolor=lightgray]{
\psline(8,3)(11,3) \psline[liftpen=0](11,2)(8,2) }

\psline[linestyle=dotted](8,5)(9,5)
\psline[linestyle=solid](9,5)(13,5)
\psline[linestyle=dotted](8,4)(9,4)
\psline[linestyle=solid](9,4)(13,4)
\psline[linestyle=dotted](8,3)(10,3)
\psline[linestyle=solid](10,3)(12,3)
\psline[linestyle=dotted](8,2)(11,2)
\psline[linestyle=solid](11,2)(12,2)

\psline[linestyle=dotted](8,2)(8,5)
\psline[linestyle=dotted](9,2)(9,4)
\psline[linestyle=solid](9,4)(9,5)
\psline[linestyle=dotted](10,2)(10,3)
\psline[linestyle=solid](10,3)(10,5)
\psline[linestyle=solid](11,2)(11,5)
\psline[linestyle=solid](12,2)(12,5)
\psline[linestyle=solid](13,4)(13,5)

\rput(10.5,1){${\bf T}_{\lambda-\mu}$}

\end{pspicture}
\qquad\quad
\begin{pspicture}(0,0)(7,4)
\psset{xunit=.5cm, yunit=.5cm}

\rput(7,6){$\lambda':=(5,5,5) \quad \quad \mu':=(1,3,4)$}

\psline[linestyle=solid](0,5)(5,5)
\psline[linestyle=solid](0,4)(5,4)
\psline[linestyle=solid](0,3)(5,3)
\psline[linestyle=solid](0,2)(5,2)

\psline[linestyle=solid](0,2)(0,5)
\psline[linestyle=solid](1,2)(1,5)
\psline[linestyle=solid](2,2)(2,5)
\psline[linestyle=solid](3,2)(3,5)
\psline[linestyle=solid](4,2)(4,5)
\psline[linestyle=solid](5,2)(5,5)

\rput(2.5,1){${\bf T}_{\lambda'}$}

\pscustom[linestyle=none,fillstyle=solid,fillcolor=lightgray]{
\psline(8,5)(9,5) \psline[liftpen=0](9,4)(8,4) }

\pscustom[linestyle=none,fillstyle=solid,fillcolor=lightgray]{
\psline(8,4)(11,4) \psline[liftpen=0](11,3)(8,3) }

\pscustom[linestyle=none,fillstyle=solid,fillcolor=lightgray]{
\psline(8,3)(12,3) \psline[liftpen=0](12,2)(8,2) }

\psline[linestyle=dotted](8,5)(9,5)
\psline[linestyle=solid](9,5)(13,5)
\psline[linestyle=dotted](8,4)(9,4)
\psline[linestyle=solid](9,4)(13,4)
\psline[linestyle=dotted](8,3)(11,3)
\psline[linestyle=solid](11,3)(13,3)
\psline[linestyle=dotted](8,2)(12,2)
\psline[linestyle=solid](12,2)(13,2)

\psline[linestyle=dotted](8,2)(8,5)
\psline[linestyle=dotted](9,2)(9,4)
\psline[linestyle=solid](9,4)(9,5)
\psline[linestyle=dotted](10,2)(10,4)
\psline[linestyle=solid](10,4)(10,5)
\psline[linestyle=dotted](11,2)(11,3)
\psline[linestyle=solid](11,3)(11,5)
\psline[linestyle=solid](12,2)(12,5)
\psline[linestyle=solid](13,2)(13,5)

\rput(10.5,1){${\bf T}_{\lambda'-\mu'}$}

\end{pspicture}
\end{center}

\noindent
By reordering the columns of ${\mathbf T}_{\lambda - \mu}$
according to their length we see that $I_{\lm}$ is isomorphic to the
Ferrers ideal associated to the partition $(\lambda_1 -\mu_1,
\ldots, \lambda_n - \mu_n)$. Note, however, that isomorphic
generalized Ferrers ideals have in general non-isomorphic
specializations.

\begin{example} \label{ex-diff-spec}
Let $\lambda := (5, 4, 4), \ \mu := (1, 2, 3)$ and $\lambda' := (5,
5, 5), \ \mu' := (1, 3, 4)$. Then the generalized Ferrers ideals
% $I_{\lm}$ and $I_{\lambda' - \mu'}$
\[
I_{\lambda-\mu}=(x_1y_2, x_1y_3, x_1y_4, x_1y_5, x_2y_3, x_2y_4,
x_3y_4),
\]
\[
I_{\lambda'-\mu'}=(x_1y_2, x_1y_3, x_1y_4, x_1y_5, x_2y_4, x_2y_5,
x_3y_5)
\]
are isomorphic, while their specializations
% $\ov{I}_{\lm}$ and $\ov{I}_{\lambda' - \mu'}$
\[
\ov{I}_{\lm} =(x_1x_2, x_1x_3, x_1x_4, x_1x_5, x_2x_3, x_2x_4,
x_3x_4),
\]
\[
\ov{I}_{\lambda' - \mu'}=(x_1x_2, x_1x_3, x_1x_4, x_1x_5,
x_2x_4, x_2x_5, x_3x_5)
\]
are not isomorphic. Indeed, the graph $\ov{G}_{\lambda' - \mu'}$
corresponding to $\ov{I}_{\lambda' - \mu'}$ has two vertices of
degree two whereas the graph $\ov{G}_{\lambda - \mu}$ corresponding
to $\ov{I}_{\lm}$ does not have any vertex of degree two.
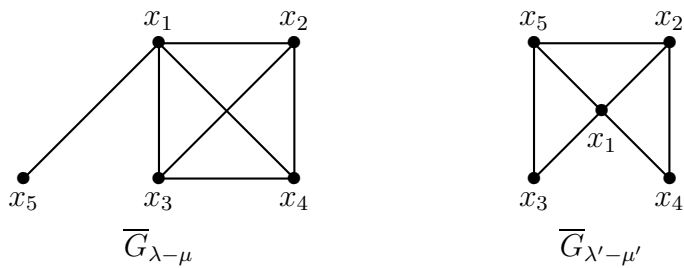
\begin{center}
\begin{pspicture}(-1,-1)(4,2.5)
\psset{xunit=.6cm, yunit=.6cm}

\rput(0,0){$\bullet$} \rput(3,0){$\bullet$} \rput(6,0){$\bullet$}
\rput(3,3){$\bullet$} \rput(6,3){$\bullet$}
\psline[linestyle=solid](0,0)(3,3)
\psline[linestyle=solid](3,0)(3,3)
\psline[linestyle=solid](3,3)(6,0)
\psline[linestyle=solid](3,0)(6,3)
\psline[linestyle=solid](3,3)(6,3)
\psline[linestyle=solid](6,3)(6,0)
\psline[linestyle=solid](3,0)(6,0)

\rput(0,-0.5){$x_5$} \rput(3,-0.5){$x_3$} \rput(6,-0.5){$x_4$}
\rput(3,3.5){$x_1$} \rput(6,3.5){$x_2$}

\rput(3,-1.5){$\overline{G}_{\lambda-\mu}$}

\end{pspicture}
\qquad\qquad
\begin{pspicture}(-1,-1)(4,2.5)
\psset{xunit=.6cm, yunit=.6cm}

\rput(0,0){$\bullet$} \rput(0,3){$\bullet$} \rput(3,0){$\bullet$}
\rput(3,3){$\bullet$}\rput(1.5,1.5){$\bullet$}

\psline[linestyle=solid](3,0)(3,3)
\psline[linestyle=solid](0,3)(3,3)
\psline[linestyle=solid](0,0)(0,3)
\psline[linestyle=solid](1.5,1.5)(0,0)
\psline[linestyle=solid](1.5,1.5)(0,3)
\psline[linestyle=solid](1.5,1.5)(3,0)
\psline[linestyle=solid](1.5,1.5)(3,3)

\rput(0,-0.5){$x_3$} \rput(0,3.5){$x_5$} \rput(3,-0.5){$x_4$}
\rput(3,3.5){$x_2$} \rput(1.5,0.8){$x_1$}

\rput(1.5,-1.5){$\overline{G}_{\lambda'-\mu'}$}

\end{pspicture}
\end{center}
\end{example}

We want to show that the specialization of a generalized Ferrers
ideal has a minimal free cellular resolution. This requires some
preparation. We begin by describing the minimal free resolution of a
generalized Ferrers ideal. The complex $X_{nm}$ has been introduced
in Example~\ref{ex-compl-bip}.

\begin{definition} \label{def-gen-cell-c}
The labeled polyhedral cell complex $X_{\lm}$ associated to
${\lambda}$ and $\mu$ is the labeled subcomplex of $X_{n,m}$
consisting of all the faces of $X_{n,m}$ whose vertices are labeled
by monomials in the generalized Ferrers ideal $I_{\lm}$.
\end{definition}

This complex captures the information about the resolution of
$I_{\lm}$. In fact, we have:

\begin{lemma} \label{lem-mfr-f}
The complex ${\mathcal F}_{X_{\lm}}$ of free $R$-modules provides
the minimal free $\ZZ^{m+n}$-graded resolution of $I_{\lm}$.
\end{lemma}

\begin{proof}
Since the generalized Ferrers ideal is isomorphic to a Ferrers
ideal, we may restrict ourselves to this case, i.e.\  $\mu = 0$, by
permuting the variables $y_1,\ldots,y_m$ suitably. For Ferrers
ideals, the claim is shown as in \cite[Theorem 3.2]{CN1}.
\end{proof}

\begin{corollary}\label{cor-genbetti}
The minimal $\ZZ$-graded free resolution of $I_{\lm}$  is $2$-linear
and, for $i > 0$, the $i${\rm{-th}} Betti number of $R/I_{\lm}$ is
given by
\[
\beta_i(R/I_{\lm}) = {\lambda_1 - \mu_1\choose i} + {\lambda_2 - \mu_2 +1
\choose i}+ \ldots + {\lambda_n - \mu_n +n-1
\choose i} - {n \choose i+1}.
\]
\end{corollary}

\begin{proof}
This follows as in \cite{CN1} because each $i$-dimensional face of
$X_{\lm}$ has a label of total degree $i+2$.
\end{proof}

Now we want to specialize. Notice that if $\mu_i \leq i-2$ for some
$i \geq 2$ and $\lambda_{i-1} \geq i$, the two monomials $x_{i-1}
y_i$ and $x_i y_{i-1}$ in $I_{\lm}$ specialize to the same monomial.
Excluding this case, we get:

\begin{lemma} \label{lem-same-gen}Suppose in addition that
$\mu_i \geq i-1$, $i= 1, \ldots, n$.  Then the ideals $I_{\lm}$ and
$\ov{I}_{\lm}$ have the same number of minimal generators, namely
$|\lambda| - |\mu| = \lambda_1 + \cdots + \lambda_n - [\mu_1 +
\cdots + \mu_n]$.
\end{lemma}

\begin{proof}
The assumption guarantees that the specialization map is injective
on the set of minimal generators of $I_{\lm}$.
\end{proof}

This observation shows that the labels of the following complex are
pairwise distinct.

\begin{definition} \label{def-spec-cell-c}
The labeled polyhedral cell complex $\ov{X}_{\lm}$ associated to
${\lambda}$ and $\mu$ is the complex obtained from $X_{\lm}$ by
specializing its labels. In particular, both complexes have the same
supporting cell complex.
\end{definition}

\begin{example}
Let $\lambda := (4, 4, 4)$ and $\mu := (1, 2, 3)$. Below we depict
the complex $X_{\lm}$ on the left-hand side and its specialization
$\ov{X}_{\lm}$ on the right-hand side.

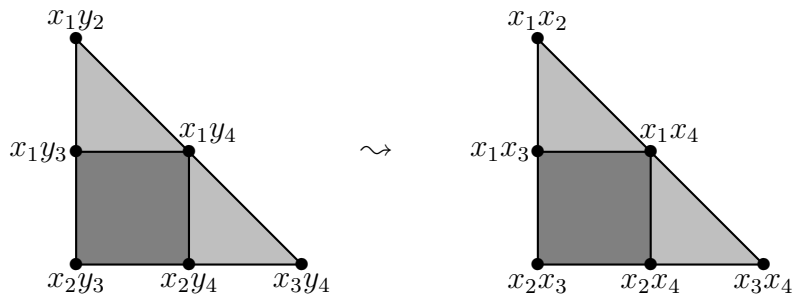
\begin{pspicture}(-3,-1)(3,4)
\psset{xunit=.5cm, yunit=.5cm}

\pscustom[linestyle=none,fillstyle=solid,fillcolor=lightgray]{
\psline(0,6)(0,3) \psline[liftpen=0](0,3)(3,3) }

\pscustom[linestyle=none,fillstyle=solid,fillcolor=lightgray]{
\psline(3,3)(3,0) \psline[liftpen=0](3,0)(6,0) }

\pscustom[linestyle=none,fillstyle=solid,fillcolor=gray]{
\psline(0,0)(3,0) \psline[liftpen=0](3,3)(0,3) }

\psline[linestyle=solid](0,0)(6,0)
\psline[linestyle=solid](0,0)(0,6)
\psline[linestyle=solid](0,6)(6,0)
\psline[linestyle=solid](0,3)(3,3)
\psline[linestyle=solid](3,3)(3,0)

\rput(0,0){$\bullet$}

\rput(3,0){$\bullet$}

\rput(6,0){$\bullet$}

\rput(0,3){$\bullet$}

\rput(3,3){$\bullet$}

\rput(0,6){$\bullet$}

\rput(3.5,3.5){$x_1y_4$}

\rput(3,-0.5){$x_2y_4$}

\rput(6,-0.5){$x_3y_4$}

\rput(0,6.5){$x_1y_2$}

\rput(-1,3){$x_1y_3$}

\rput(0,-0.5){$x_2y_3$}

\rput(8,3){$\leadsto$}
\end{pspicture}
\begin{pspicture}(-3,-1)(3,4)
\psset{xunit=.5cm, yunit=.5cm}

\pscustom[linestyle=none,fillstyle=solid,fillcolor=lightgray]{
\psline(0,6)(0,3) \psline[liftpen=0](0,3)(3,3) }

\pscustom[linestyle=none,fillstyle=solid,fillcolor=lightgray]{
\psline(3,3)(3,0) \psline[liftpen=0](3,0)(6,0) }

\pscustom[linestyle=none,fillstyle=solid,fillcolor=gray]{
\psline(0,0)(3,0) \psline[liftpen=0](3,3)(0,3) }

\psline[linestyle=solid](0,0)(6,0)
\psline[linestyle=solid](0,0)(0,6)
\psline[linestyle=solid](0,6)(6,0)
\psline[linestyle=solid](0,3)(3,3)
\psline[linestyle=solid](3,3)(3,0)

\rput(0,0){$\bullet$}

\rput(3,0){$\bullet$}

\rput(6,0){$\bullet$}

\rput(0,3){$\bullet$}

\rput(3,3){$\bullet$}

\rput(0,6){$\bullet$}

\rput(3.5,3.5){$x_1x_4$}

\rput(3,-0.5){$x_2x_4$}

\rput(6,-0.5){$x_3x_4$}

\rput(0,6.5){$x_1x_2$}

\rput(-1,3){$x_1x_3$}

\rput(0,-0.5){$x_2x_3$}

\end{pspicture}

\noindent The facets of both complexes are two triangles and one
rectangle.
\end{example}

The main result of this note is:

\begin{theorem} \label{thm-mfr-spec}
If  $\mu_i \geq i-1$ $(i = 1,\ldots,n)$, then the complex ${\mathcal
F}_{\ov{X}_{\lm}}$ of free $S$-modules provides the minimal free
$\ZZ^m$-graded resolution of the specialization of the generalized
Ferrers ideal $\ov{I}_{\lm}$.
\end{theorem}

\begin{proof}
Our strategy is to reduce the claim to the corresponding statement
for generalized Ferrers ideeals by applying the  criterion of Bayer
and Sturmfels (see Lemma \ref{lem-acyc-crit}) twice.

Let $\ov{{\mathbf c}} := (c_1,\ldots,c_n) \in \NN_0^m$ and denote by
$\ov{{\mathbf c'}} \in \ZZ^m$ the degree of the least common
multiple of the labels of the vertices in
$(\ov{X}_{\lm})_{\preccurlyeq \ov{{\mathbf c}}}$. Then we get
$(\ov{X}_{\lm})_{\preccurlyeq \ov{{\mathbf c}}} =
(\ov{X}_{\lm})_{\preccurlyeq \ov{{\mathbf c'}}}$. Now we define
${\mathbf c} := (a_1,\ldots,a_n, b_1,\ldots,b_m)$ by
\[
a_i = \max \{0, c_i - 1\},  \quad i = 1,\ldots,n,
\]
and
\[
b_i = \left \{\begin{array}{ll}
c_i & \mif \ m < i \leq n \\
c_i - a_i & \mif \ 1 \leq i \leq n.
\end{array}
\right.
\]
The crucial observation is:

\begin{quote}
{\it Claim:} $(\ov{X}_{\lm})_{\preccurlyeq \ov{{\mathbf c'}}}$ is
the labeled cell complex obtained from $(X_{\lm})_{\preccurlyeq
{\mathbf c}}$ \newline by specializing its labels.
\end{quote}

\noindent Indeed,  observe that the specialization of the monomial
$x^{{\mathbf a}} y^{{\mathbf b}} \in R$ is the monomial $x^{{\mathbf
a'} + {\mathbf b}} \in S$, where ${\mathbf a'} \in \ZZ^m$ is the
vector obtained from ${\mathbf a} \in \ZZ^n$ by appending it with
$n-m$ zero entries. This provides the claim.

According to Lemma \ref{lem-mfr-f}, the complex ${\mathcal
F}_{X_{\lm}}$ is exact. Thus Lemma \ref{lem-acyc-crit} yields that
the complex $(X_{\lm})_{\preccurlyeq {\mathbf c}}$ is acyclic over
$K$. Hence, the above claim shows that $(\ov{X}_{\lm})_{\preccurlyeq
\ov{{\mathbf c}}} = (\ov{X}_{\lm})_{\preccurlyeq \ov{{\mathbf c'}}}$
is also acyclic. Applying the Bayer-Sturmfels criterion, Lemma
\ref{lem-acyc-crit}, now to ${\mathcal F}_{\ov{X}_{\lm}}$ completes
the proof.
\end{proof}

\begin{corollary}
If  $\mu_i \geq i-1$ $(i = 1, \ldots, n)$, then the ideal
$\ov{I}_{\lm}$ has a 2-linear $\ZZ$-graded free resolution, i.e. its
Castelnuovo-Mumford regularity is two.
\end{corollary}

\begin{proof}
This follows by combining Corollary \ref{cor-genbetti} and the claim
in the above proof.
\end{proof}

Using \cite{EGHP}, Proposition 0.3, the last result implies in
particular that  $I_{\lm}$ is the homogeneous ideal of a small
subscheme in $\PP^{m-1}$ that is not necessarily reduced. In
\cite{EGHP}, Theorem 6.1,  Eisenbud, Green, Hulek, and Popescu
construct a free resolution for every reduced subscheme that is the
union of linear subspaces and that has a $2$-linear free resolution.
However, in general this resolution is not minimal though it gives
the exact number of minimal generators of the homogeneous ideal
$I_X$. Our Theorem~\ref{thm-mfr-spec} treats  cases where $I_X$ is a
not necessarily reduced monomial ideal, and it has a stronger
conclusion.

We can also interpret the above results using the concept of
lifting. Indeed, let $I$ be an ideal in the commutative ring $A$ and
let $u_1,\ldots,u_t$ be elements in $A$. Set $B :=
A/(u_1,\ldots,u_t)A$ and let $J \subset A$ be an ideal. Then $I$ is
said to be a {\it $t$-lifting} of $J$ if $\{u_1,\ldots,u_t\}$ is an
$A/I$-regular sequence and $(I, u_1,\ldots,u_t)/(u_1,\ldots,u_t)
\cong J$ (see \cite{MN-lifting}, Definitions 2.1 and 2.3).

Recall our assumption $m \geq n$. Hence $R$ is a subring of the
polynomial ring $R' := K[x_0,\ldots,x_m, y_0,\ldots,y_m]$.

\begin{corollary} \label{cor-lifting}
If $\mu_i \geq i-1$ $(i = 1,\ldots,n)$, then the ideal $I_{\lm} R'$
is an $m$-lifting of the ideal $\ov{I}_{\lm} \subset S$.
\end{corollary}

\begin{proof}
Obviously, we have $(I_{\lm} R' + (y_1 - x_1,\ldots,y_m - x_m))/(y_1
- x_1,\ldots,y_m - x_m) \cong \ov{I}_{\lm}$.  It remains to show
that $\{y_1 - x_1,\ldots,y_m - x_m\}$ is an $R'/I_{\lm} R'$- regular
sequence. The minimal free resolution of $I_{\lm}$ over $R$ has the
same length as the minimal free resolution of $I_{\lm} R'$ over
$R'$. Moreover, Theorem \ref{thm-mfr-spec} shows that both
resolutions have the same length as the minimal free resolution of
$\ov{I}_{\lm}$ over $S$. Hence, the Auslander-Buchsbaum formula
provides that
\[
\depth R'/I_{\lm} R' = m + \depth S/\ov{I}_{\lm} = m + \depth
R'/I_{\lm} + (y_1 - x_1, \ldots, y_m - x_m).
\]
It follows that $\{y_1 - x_1,\ldots,y_m - x_m\}$ is an $R'/I_{\lm}
R'$- regular sequence.
\end{proof}

Probably, the last result could be shown directly by brute force,
thus giving an alternative approach to the results about the
resolutions of the specializations. However, our above approach
seems more elegant and transparent.

\section{Threshold graphs and stable ideals} \label{sec-thres}

We are going to discuss the graphs and ideals, respectively, that we
obtain as specializations of Ferrers graphs and ideals. Allowing
loops, each graph $G$ on the vertex set $[m] = \{1,\ldots,m\}$
defines the edge ideal $I_G \subset S$ that is generated by the
monomials $x_i x_j$ such that $(i, j)$ is an edge of $G$. This
provides a one-to-one correspondence between graphs on $[m]$ and
monomial ideals in $S$ whose minimal generators all have degree two.

Consider now a graph $G$ on $[m]$ without isolated vertices. This
assumption is harmless as far as the edge ideal is concerned. Order
the vertices of $G$ as follows. Denote by 1 one of the vertices of
highest degree. Assume we have chosen vertices $1,\ldots,i-1$ where
$2 \leq i \leq m$. Then we denote by $i$ one of the vertices of
highest degree of the subgraph of $G$ on the vertex set
$\{i,\ldots,m\}$. Now we define
\[
n := \max \{i \s \mbox{There is a vertex  $j \geq i$ such that $(i,
j)$ is an edge of G}\}.
\]
Furthermore, we set $\lambda := (\lambda_1,\ldots,\lambda_n)$ where
\[
\lambda_i := \max \{j \s \mbox{$(i, j)$ is an edge of G}\}.
\]
Note that $\lambda_n \geq n$ by the definition of $n$. Assume that
$\lambda_1 = m$ and that  $\lambda_1 \geq \lambda_2 \geq \cdots \geq
\lambda_n$. Then we can define $\mu := (\mu_1,\ldots,\mu_n)$ where
\[
\mu_i := -1 + \min \{j \geq i \s \mbox{$(i, j)$ is an edge of G}\}.
\]

\begin{example} \label{ex-al-compl}
Let $G$ be the graph obtained from the complete graph on 4 vertices
by taking away one edge. Then the above procedure gives
$I_G=I_{\lambda-\mu}$, where $\lambda = (4, 4),\; \mu = (1, 2)$.
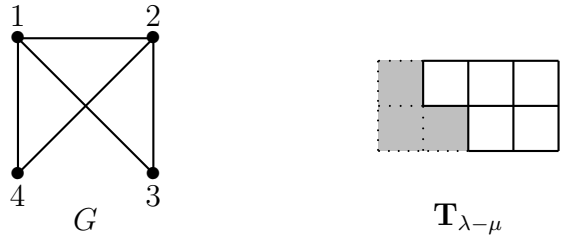
\begin{center}
\begin{pspicture}(0,-1)(4,2.5)
\psset{xunit=.6cm, yunit=.6cm}

\rput(3,0){$\bullet$} \rput(6,0){$\bullet$} \rput(3,3){$\bullet$}
\rput(6,3){$\bullet$} \psline[linestyle=solid](3,0)(3,3)
\psline[linestyle=solid](3,3)(6,0)
\psline[linestyle=solid](3,3)(6,3)
\psline[linestyle=solid](6,3)(3,0)
\psline[linestyle=solid](6,3)(6,0)

\rput(3,-0.5){$4$} \rput(6,-0.5){$3$} \rput(3,3.5){$1$}
\rput(6,3.5){$2$}

\rput(4.5,-1){$G$}
\end{pspicture}
\qquad\qquad
\begin{pspicture}(4,-1)(8,2.5)
\psset{xunit=.6cm, yunit=.6cm}

\pscustom[linestyle=none,fillstyle=solid,fillcolor=lightgray]{
\psline(8,2.5)(9,2.5) \psline[liftpen=0](9,1.5)(8,1.5) }

\pscustom[linestyle=none,fillstyle=solid,fillcolor=lightgray]{
\psline(8,1.5)(10,1.5) \psline[liftpen=0](10,0.5)(8,0.5) }

\psline[linestyle=dotted](8,2.5)(9,2.5)
\psline[linestyle=solid](9,2.5)(12,2.5)
\psline[linestyle=dotted](8,1.5)(9,1.5)
\psline[linestyle=solid](9,1.5)(12,1.5)
\psline[linestyle=dotted](8,0.5)(10,0.5)
\psline[linestyle=solid](10,0.5)(12,0.5)

\psline[linestyle=dotted](8,0.5)(8,2.5)
\psline[linestyle=dotted](9,0.5)(9,1.5)
\psline[linestyle=solid](9,1.5)(9,2.5)
\psline[linestyle=solid](10,0.5)(10,2.5)
\psline[linestyle=solid](11,0.5)(11,2.5)
\psline[linestyle=solid](12,0.5)(12,2.5)

\rput(10,-1){${\bf T}_{\lambda-\mu}$}

\end{pspicture}
\end{center}
\end{example}

\begin{proposition} \label{prop-graphs}
Adopt the above notation and assume that the graph $G$ satisfies
\[
\mu_1 \geq \mu_2 \geq \cdots \geq \mu_n
\]
and
\begin{equation} \label{cond}
(i, j) \mbox{ is an edge of $G$ whenever } 1 \leq i \leq n \
\mbox{and } \mu_i < j \leq \lambda_i.
\end{equation}
Then its edge ideal is $I_G = \ov{I}_{\lm}$ and ${\mathcal
F}_{\ov{X}_{\lm}}$ is the  minimal free cellular $\ZZ^{m}$-graded
resolution of $I_G$. In particular,
\begin{eqnarray*}
\reg (I_G) & = & 2, \\
\htt I_G  & = & \min \{ \min_j \{ \lambda_j - \mu_j +j-1 \}, n \}, \\
\depth S/I_G  & = & m - \max_j\{ \lambda_j - \mu_j +j-1 \}
\end{eqnarray*}
and the $i${\rm{-th}} Betti number of $S/I_{G}$ is given by
\[
\beta_i(S/I_{G}) = {\lambda_1 - \mu_1\choose i} + {\lambda_2 - \mu_2 +1
\choose i}+ \ldots + {\lambda_n - \mu_n +n-1
\choose i} - {n \choose i+1}.
\]
\end{proposition}

\begin{proof}
This follows by Theorem \ref{thm-mfr-spec} and Corollary
\ref{cor-lifting} from the corresponding results for Ferrers ideals
in \cite{CN1}.
\end{proof}

\begin{example} \label{ex-exception}
Let $G'$ be the graph obtained from the complete graph on 4 vertices
by taking away two edges that share a common vertex. Then the above
procedure shows that Proposition \ref{prop-graphs} applies to $G'$
with  $n = 2$ and $\lambda = (4, 3),\; \mu = (1, 2)$.
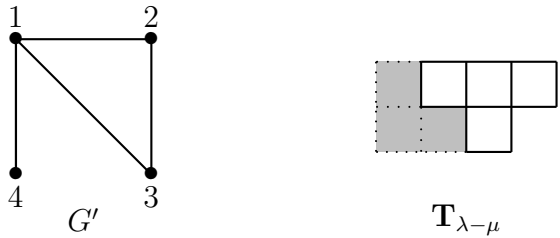
\begin{center}
\begin{pspicture}(0,-1)(4,2.5)
\psset{xunit=.6cm, yunit=.6cm}

\rput(3,0){$\bullet$} \rput(6,0){$\bullet$} \rput(3,3){$\bullet$}
\rput(6,3){$\bullet$} \psline[linestyle=solid](3,0)(3,3)
\psline[linestyle=solid](3,3)(6,0)
\psline[linestyle=solid](3,3)(6,3)
\psline[linestyle=solid](6,3)(6,0)
\psline[linestyle=solid](3,0)(3,3)

\rput(3,-0.5){$4$} \rput(6,-0.5){$3$} \rput(3,3.5){$1$}
\rput(6,3.5){$2$}

\rput(4.5,-1){$G'$}
\end{pspicture}
\qquad\qquad
\begin{pspicture}(4,-1)(8,2.5)
\psset{xunit=.6cm, yunit=.6cm}

\pscustom[linestyle=none,fillstyle=solid,fillcolor=lightgray]{
\psline(8,2.5)(9,2.5) \psline[liftpen=0](9,1.5)(8,1.5) }

\pscustom[linestyle=none,fillstyle=solid,fillcolor=lightgray]{
\psline(8,1.5)(10,1.5) \psline[liftpen=0](10,0.5)(8,0.5) }

\psline[linestyle=dotted](8,2.5)(9,2.5)
\psline[linestyle=solid](9,2.5)(12,2.5)
\psline[linestyle=dotted](8,1.5)(9,1.5)
\psline[linestyle=solid](9,1.5)(12,1.5)
\psline[linestyle=dotted](8,0.5)(10,0.5)
\psline[linestyle=solid](10,0.5)(11,0.5)

\psline[linestyle=dotted](8,0.5)(8,2.5)
\psline[linestyle=dotted](9,0.5)(9,1.5)
\psline[linestyle=solid](9,1.5)(9,2.5)
\psline[linestyle=solid](10,0.5)(10,2.5)
\psline[linestyle=solid](11,0.5)(11,2.5)
\psline[linestyle=solid](12,1.5)(12,2.5)

\rput(10,-1){${\bf T}_{\lambda-\mu}$}

\end{pspicture}
\end{center}
Notice that $G'$ with a particular labeling is the graph that is excluded as a pattern of the
graphs considered in \cite{H}.
\end{example}

\begin{corollary} \label{cor-CM}
Adopt the notation and assumptions of {\rm Proposition
\ref{prop-graphs}}. Then $S/I_G$ is a Cohen-Macaulay ring if and
only if
\[
\min_j \{\lambda_j -\mu_j + j - 1 \}  = \max _j \{\lambda_j -\mu_j +
j - 1 \} \leq n.
\]
\end{corollary}

\begin{proof}
$S/I_G$ is Cohen-Macaulay if and only $\dim S/I_G = \depth S/I_G$.
Hence Proposition \ref{prop-graphs} provides the claim.
\end{proof}

We now discuss classes of graphs or ideals to which
Proposition~\ref{prop-graphs} applies. Recall that a monomial ideal
$I\subset S$ is called {\it strongly stable} if $x_i
\displaystyle\frac{x^{\mathbf a}}{x_j} \in I$ whenever $x^{\mathbf
a} \in I$, $x_j$ divides $x^{\mathbf a}$, and $1 \leq i < j$. The
squarefree monomial ideal $I$ is said to be {\it squarefree strongly
stable} if  $x_i \displaystyle\frac{x^{\mathbf a}}{x_j} \in I$
whenever $x^{\mathbf a}$ is a minimal generator of $I$, $x_j$
divides $x^{\mathbf a}$, $x_i$ does not divide $x^{\mathbf a}$, and
$1 \leq i < j$. Note that there are also the, in general,  weaker
conditions of being stable or squarefree stable. However, for ideals
generated in degree two, the corresponding concepts are equivalent.
Eliahou and Kervaire describe  in \cite{EK} the minimal free
resolution of an arbitrary strongly stable ideal $I$. If $I$ is
generated in degree two, our results show that $I$ admits a cellular
resolution. More precisely:

\begin{example} \label{ex-stable}
Let $I \subset S$ be a strongly stable ideal whose minimal
generators have degree two and such that $x_1 x_m \in I$. Let $G$ be
the corresponding graph. Then the stability property guarantees that
$G$ satisfies Condition (\ref{cond}) where $\mu_i = i-1$ for
$i=1,\ldots,n$. In particular, Corollary \ref{cor-CM} immediately
implies the well-known fact that $S/I$ is Cohen-Macaulay if and only
if $m = \lambda_1 = \cdots \lambda_n = n$, that is $I =
(x_1,\ldots,x_m)^2$.
\end{example}

Recall that a graph $G$ on $[m]$ is called a {\em threshold graph}
if there is a vector $w = (w_1,\ldots,w_m) \in \RR^m$ such that $(i,
j)$ is an edge of $G$ if and only if $w_i + w_j > 0$. We refer to
the book by Mahadev and Peled \cite{MP} for a wealth of information
on threshold graphs and to the work of Klivans and Reiner \cite{KR}
for many alternative characterizations of threshold graphs.

\begin{corollary} \label{cor-thres}
Let $G$ be a threshold graph on $[m]$ and denote by $\lambda_i$ the
degree of the vertex $i$. Order the vertices such that $\lambda_1
\geq \lambda_2 \geq \cdots \geq \lambda_m \geq 1$. Define
$$
n := \max \{i \s \lambda_i \geq i+1\}
$$
and
$$
\mu = (\mu_1,\ldots,\mu_n) := (1, 2,\ldots,n).
$$
Then the edge ideal of $G$  is $I_G = \ov{I}_{\lm}$ and ${\mathcal
F}_{\ov{X}_{\lm}}$ is the  minimal free cellular $\ZZ^{m}$-graded
resolution of $I_G$. In particular, $\reg (I_G)  =  2, \; \htt I_G =
n, \; \depth S/I_G   = 1$ and the $i${\rm{-th}} Betti number of
$S/I_{G}$ is given by
\[
\beta_i(S/I_{G}) = {\lambda_1 - 1\choose i} + {\lambda_2 - 1
\choose i}+ \ldots + {\lambda_n - 1
\choose i} - {n \choose i+1}.
\]
\end{corollary}

\begin{proof}
Since $G$ is threshold, it satisfies the conditions in Proposition
\ref{prop-graphs}. In fact, the tableaux ${\mathbf T}_{\lm}$
corresponds to the so-called up-degree sequence of $G$.
\end{proof}

\begin{remark}
(i) Corollary \ref{cor-thres} shows explicitly how each threshold
graph can be obtained as the specialization of a Ferrers graphs.
Thus, it explains the similar algebraic properties of Ferrers and
threshold graphs.

(ii) It is known (see \cite{KR}) that a graph without loops is
threshold if and only if it is shifted. Equivalently, the edge
ideals of threshold graphs are precisely the squarefree strongly
stable ideals that are generated in degree two. The minimal free
resolution of an arbitrary squarefree monomial ideal has been
described by Aramova, Herzog and Hibi in \cite{AHH}.
\end{remark}

We end our note by remarking that, by suitably modifying the vertex
labels if necessary, our methods apply to more graphs than discussed
so far.

\begin{example}
Consider the graph $\ov{G}_{\lambda' - \mu'}$ that is described in
Example~\ref{ex-diff-spec}. This is the same graph discussed in
\cite[Example 4.3]{H}, but with a different labeling. As remarked
earlier, its edge ideal is the specialization of $I_{\lambda'-\mu'}$
with $\lambda' = (5, 5, 5)$ and $\mu' = (1, 3, 4)$.  Hence the
cellular resolution of $\ov{G}_{\lambda' - \mu'}$ is given by the
polyhedral cell complex pictured below
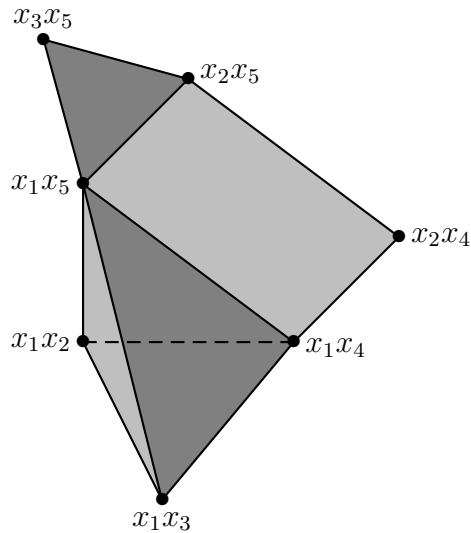
\begin{center}
\begin{center}
\begin{pspicture}(-1,-2.5)(5,4.5)
\psset{xunit=.7cm, yunit=.7cm}

\pscustom[linestyle=none,fillstyle=solid,fillcolor=lightgray]{
\psline(0,0)(0,3) \psline[liftpen=0](0,3)(1.5,-3) }

\pscustom[linestyle=none,fillstyle=solid,fillcolor=gray]{
\psline(1.5,-3)(0,3) \psline[liftpen=0](0,3)(4,0) }

\pscustom[linestyle=none,fillstyle=solid,fillcolor=lightgray]{
\psline(0,3)(2,5) \psline[liftpen=0](6,2)(4,0) }

\pscustom[linestyle=none,fillstyle=solid,fillcolor=gray]{
\psline(0,3)(-0.75,5.75) \psline[liftpen=0](-0.75,5.75)(2,5) }

\psline[linestyle=dashed](0,0)(4,0)
\psline[linestyle=solid](0,0)(0,3)
\psline[linestyle=solid](0,3)(4,0)
\psline[linestyle=solid](0,0)(1.5,-3)
\psline[linestyle=solid](1.5,-3)(4,0)
\psline[linestyle=solid](1.5,-3)(0,3)

\psline[linestyle=solid](4,0)(6,2)
\psline[linestyle=solid](0,3)(2,5)
\psline[linestyle=solid](6,2)(2,5)
\psline[linestyle=solid](0,3)(-0.75,5.75)
\psline[linestyle=solid](2,5)(-0.75,5.75)

\rput(0,0){$\bullet$} \rput(0,3){$\bullet$} \rput(4,0){$\bullet$}
\rput(1.5,-3){$\bullet$} \rput(6,2){$\bullet$} \rput(2,5){$\bullet$}
\rput(-0.75,5.75){$\bullet$}

\rput(-0.75,6.15){$x_3x_5$} \rput(1.5,-3.4){$x_1x_3$}
\rput(-0.8,0){$x_1x_2$} \rput(-0.8,3){$x_1x_5$}
\rput(4.8,-0.1){$x_1x_4$} \rput(6.8,2){$x_2x_4$}
\rput(2.8,5.1){$x_2x_5$}

\end{pspicture}
\end{center}
\end{center}
This cell complex has a $2$-simplex as a facet, whereas Horwitz's
cell complex does not have such a facet. This shows in particular
that the maps in the free resolutions constructed by Horwitz and by
our methods are in general different.

After completing the first version of this note, Horwitz  informed us that by applying his methods to the graph $\ov{G}_{\lambda' - \mu'}$ with the above labeling, he gets the same abstract cell complex, but with a different labeling of its vertices, so again his resulting maps in the free resolution are different from ours.
\end{example}

\medskip

\noindent {\bf Acknowledgement.} The authors would like to thank Vic
Reiner for inspiring discussions.

\end{document}